\documentclass[11pt,fullpage]{article}

\usepackage[pdftex]{graphicx}
\usepackage{latexsym}
\usepackage{color}
\usepackage{multirow}
\usepackage{amsmath}
\usepackage{amssymb}
\usepackage{amsthm}
\usepackage{amsfonts}
\usepackage{bm}

\newtheorem{theorem}{Theorem}

\newtheorem{assumption}{Assumption}
\newtheorem{corollary}{Corollary}

\long\def\ignore#1{}

\newcommand{\be}{\begin{equation}}
\newcommand{\ee}{\end{equation}}
\newcommand{\beqn}{\begin{eqnarray}}
\newcommand{\eeqn}{\end{eqnarray}}

\newcommand{\bfm}[1]{\mbox{\boldmath{$#1$}}}

\newcommand{\btheta}{\bfm{\theta}}

\newcommand{\by}{{\bf y}}

\newcommand{\bx}{{\bf x}}
\newcommand{\bX}{{\bf X}}

\newcommand{\bu}{{\bf u}}
\newcommand{\bj}{{\bf j}}

\newcommand{\mE}{\mathbb{E}}
\newcommand{\mP}{\mathbb{P}}

\newcommand{\cG}{{\cal G}}

\newcommand{\cJ}{{\cal J}}

\newcommand{\cE}{{\cal E}}
\newcommand{\cF}{{\cal F}}

\newcommand{\cX}{{\cal X}}

\newcommand{\cY}{{\cal Y}}

\newcommand{\sparsity}{\Delta}
\newcommand{\complexity}{r}

\newcommand{\fs}{\lfloor s \rfloor}

\long\def\ignore#1{}

\begin{document}
\title{\bf Statistical learning by sparse deep neural networks}
\author{{\bf Felix Abramovich} \\
Department of Statistics and Operations Research \\
Tel Aviv University \\
Israel\\
felix@tauex.tau.ac.il
}

\date{}

\maketitle

\begin{abstract}
We consider a deep neural network estimator based on empirical risk minimization with $l_1$-regularization. We derive a general bound for its excess risk in regression and classification (including multiclass), and prove that it is adaptively  nearly-minimax (up to log-factors) simultaneously across the entire range of various function classes.
\end{abstract}

\section{Introduction}
Deep neural networks (DNNs) have become one of the main tools in various machine learning and artificial intelligence applications, including computer vision, speech recognition, natural language processing, and others. Despite substantial efforts and contributions to the algorithmic and computational facets of deep learning, its comprehensive rigorous mathematical and statistical theory remains incomplete. We can refer to Bartlett {\em et al.} (2021) and G\"uhring {\em et al.} (2022) with references therein for reviews of the state-of-the art. This paper represents another step towards filling this gap.

Consider a statistical learning setup, where $(\bX, Y)$ have a joint distribution $\mP$ on $\cX \times \cY$  and we want
to find a prediction function $g: \cX \rightarrow \cY$ w.r.t. the loss $\ell(g(\bx),y)$. The two key examples are regression, where $\cY \subseteq \mathbb{R}^1$, $g:\cX \rightarrow \mathbb{R}^1$ with the quadratic loss $\ell(g(\bx),y)=(g(\bx)-y)^2$, and classification, where $\cY=\{0,1\}$ (binary classification) or more generally, $\cY=\{1,\ldots,K\}$ (multiclass classification), a classifier $g(\bx): \cX \rightarrow \{1,\ldots,K\}$ and the $0/1$-loss $\ell(g(\bx),y)=I\{g(\bx) \neq y\}$. 

The accuracy of $g$ is measured by its risk $R(g)=\mE\ell(g(\bX),Y)$, where the expectation is w.r.t. $\mP$, and 
the optimal predictor is $g^*=\inf_g R(g)$. 
In particular, for regression, $g^*(\bx)=\mE(Y|(\bX=\bx))$, while for classification,
$g^*(\bx)=\arg \max_{1 \leq k \leq K} P(Y=k|(\bX=\bx))$.

The conditional distribution of $Y|\bX$ is usually, however, not fully known and one should 
design a predictor $\widehat{g}$ based on the available data $D$: a random sample $(\bX_1,Y_1),\ldots,(\bX_n,Y_n)$ from $\mP$. The corresponding (conditional) risk is 
$R(\widehat{g})=\mE(\ell(\widehat{g}(\bX),Y)|D)$
and the goodness of $\widehat{g}$ is evaluated by the excess risk $\cE(\widehat{g},g^*)=\mE R(\widehat{g})-R(g^*)$.

It might be too ambitious to consider all possible predictors, so typically, one searches for a predictor within some family $\cF$. The excess risk $\cE(\widehat{g},g^*)$ can then be decomposed as
\be \label{eq:decomp}
\cE(\widehat{g},g^*)=\mE R(\widehat{g})-R(g^*)=
\left(\mE R(\widehat{g})-R(g^*_\cF)\right) + \left(R(g^*_\cF)-R(g^*)\right),
\ee
where $g^*_\cF=\arg \min_{g \in \cF} R(g)$ is the oracle choice for $g$ within $\cF$. The first term in the RHS of (\ref{eq:decomp}) represents the estimation error, whereas the second term corresponds to the approximation error reflecting the ability of functions in $\cF$ to approximate the true $g$. Enriching $\cF$ may improve the approximation error but will increase the estimation error. The choice of a proper family of estimators $\cF$ is important and may depend on a particular problem at hand.
In this paper we consider DNN estimators.

A feed-forward fully-connected DNN with $L$ layers is of the form
\begin{equation} \label{eq:dnn}
g_{\Theta}(\bx)=W_L \sigma(W_{L-1} \ldots \sigma(W_0 \bx)),
\end{equation}
where $\sigma(\cdot)$ is a (nonlinear) activation function applied element-wise, 
$W_l \in \mathbb{R}^{d_l \times d_{l-1}},\;l=1,\ldots,L$ with $d_0=d$
and $d_L=dim(\cY)$ are weight matrices and the set of parameters (weights) $\Theta=(W_0,\ldots,W_L)$. In this paper we consider the ReLU activation function $\sigma(x)=x_+$. For notational simplicity we omit bias vectors that can be incorporated into weight matrices. The total number of parameters of the DNN (\ref{eq:dnn}) is $S=\sum_{l=1}^L d_ld_{l-1}$. 

Approximate $g$  by
$g_\Theta$ from the set of DNNs $\cF(L,S)=\{\Theta: L(\Theta)=L, S(\Theta)=S)$ with a given 
architecture and estimate the unknown parameters from the training sample $(\bX_1,Y_1),\ldots,(\bX_n,Y_n)$ by empirical risk minimization (ERM) 
$$
\widehat{\Theta}=\arg\min_{\Theta \in \cF(L,S)} \frac{1}{n} \sum_{i=1}^n \ell\left(Y_i,g_\Theta(\bX_i)\right),
$$
where the loss $\ell(y,g_\Theta(\bx))$ is typically quadratic for regression and some convex surrogate for $0/1$-loss for classification.

However, the number of parameters $S$ in a DNN is typically very large leading to a severe ``curse of dimensionality'' problem and resulting in overfitting and poor generalization errors. Nevertheless, owing to their high expressive power, DNNs enable {\em sparse} approximations for various types of functions.
See G\"uhring {\em et al.} (2022) for a comprehensive review with references therein and Section \ref{subsec:adaptive} below for more references and details. 

One can think about several possible notions of DNN sparsity.
In this paper we consider the {\em connection} sparsity, where it is assumed that most of the connections between notes (weights) in $\Theta$ are not active and the number of nonzero parameters $||\Theta||_0 \ll S$. Other related types of DNN sparsity are discussed in Wen {\em et al.} (2016) and Lederer (2023). See also more remarks in Section \ref{subsec:sparsity} below.

Assume that $g$ belongs to a certain class of functions $\cG$ and
define the set of sparse DNNs $\cF(L_\cG,S_\cG)$ w.r.t. to $\cG$. Find the ERM DNN estimator within $\cF(L_\cG,S_\cG)$:
\begin{equation} \label{eq:notadaptive}
\widehat{\Theta}=\arg\min_{\Theta \in \cF(L_\cG,S_\cG)}\frac{1}{n} \sum_{i=1}^n \ell\left(Y_i,g_\Theta(\bX_i)\right)
\end{equation}

It is known that $\widehat{g}_{\widehat{\Theta}}$  in (\ref{eq:notadaptive}) with the properly chosen architecture $\cF(L_\cG,S_\cG)$  is a minimax estimator for various 
function classes, e.g. H\'older $C^s$ and Sobolev $W^{s,p}$ classes of smooth functions with regularity $s$ (Yarotsky, 2017; Schmidt-Hieber, 2020; Kim {\em et al.}, 2021), more general Besov classes $B^s_{p,q}$ (Suzuki, 2019), composition structured functions (Schmidt-Hieber, 2020; Kohler and Langer, 2021) and a certain class of piecewise smooth functions (see Imaizumi and Fukumizu, 2022 and Section \ref{subsubsec:piecewise} below for details). 

The estimator $\widehat{g}_{\widehat{\Theta}}$ in (\ref{eq:notadaptive}) is, however, not adaptive to the (usually unknown) function class and assumes that  $\cF(L_\cG,S_\cG)$ is set in advance w.r.t. a particular $\cG$.
Our objective is to design a  DNN estimator with adaptively adjusted architecture that attains the nearly-minimax risk (up to possible extra log-factors) across the entire range of various $\cG$. We show that the {\em regularized} ERM DNN estimator with Lasso $l_1$-regularization achieves this goal.

The paper is organized as follows. In Section \ref{sec:regularization} we present a $l_1$-regularized DNN estimator. We establish its near-minimaxity across various function classes for regression in Section \ref{sec:regression} and classification in Section \ref{sec:classification}. Some concluding remarks and possible future challenges are discussed in Section \ref{sec:remarks}. All proofs are given in Appendix.

\section{$l_1$-regularized DNN learning} \label{sec:regularization}
We begin with introducing some notations used throughout the paper.
Let $||\btheta||_p,\;p\geq 0$ be the $l_p$-norm of a vector $\btheta$, where, as usual, the quasi-norm $||\btheta||_0$ is the number of nonzero entries of $\btheta$. For a given DNN $\Theta \in \cF(L,S)$ it will be sometimes more convenient to re-write it in terms of a single vector $\btheta \in \mathbb{R}^S$ and denote $||\Theta||_p=||\btheta||_p$. 
Let $||g||_{L_2} = \sqrt{\int g^2(\bx)d\bx}$
be a standard $L_2$-norm of a function $g$ and $||g||_{L_2(\mP_X)} = 
\sqrt{\int g^2(\bx)d\mP_X}$
be the $L_2$-norm of g weighted w.r.t. the marginal distribution $\mP_X$ of $\bX$. For a vector-valued function $g(\bx)=(g_1(\bx),\ldots,g_K(\bx))^T$,
$||g||_{L_2}=\sqrt{\sum_{k=1}^K ||g_k||^2_{L_2}}$ and $||g||_{L_2(\mP_X)}=\sqrt{\sum_{k=1}^K ||g_k||^2_{L_2(\mP_X)}}$~.
Recall also that we use the ReLU activation function in this paper.

Consider the class $\cF(\ln n, n)$ of weakly deep DNNs with $\ln n$ layers and $n$ weights.
Define the sparse-penalized ERM DNN estimator with $l_1$-regularization (Lasso) penalty 
\begin{equation} \label{eq:lasso}
\widehat{\Theta}=\arg\min_{\Theta \in \cF(\ln n,n)}\frac{1}{n} \sum_{i=1}^n \ell\left(Y_i,g_\Theta(\bX_i)\right) + \lambda ||\Theta||_1,
\end{equation}
where $\lambda$ is a tuning parameter. The $l_1$-penalty induces connection sparsity of the resulting  $\widehat{\Theta}$ and is a widely used regularization techniques to avoid overfitting in DNN models.

While Lasso has been extensively studied for various linear and GLM setups (B\"uhlmann and van de Geer, 2011), 
its theoretical foundations within DNN framework have not been thoroughly explored yet. We can mention the recent works of Taheri {\em at el.} (2021) and Lederer (2023) although they focused on $l_1$-types constraint optimization rather than regularization.

To investigate the Lasso DNN estimator
(\ref{eq:lasso}) we need several assumptions given below.

\begin{assumption} \label{as:px} Assume that for all (generally dependent) components $X_j$'s of a random features vector $\bX \in \mathbb{R}^d$,
\begin{enumerate}
\item $\mP_X$ has a density function $f_X$ which is bounded almost surely: $f_X(\bx) \leq C\;a.s.$
\item $\mE_X X_j^2=1$ ($X_j$'s are scaled)
\item there exist constants $\kappa_1, \kappa_2, w > 1$ and $\gamma \geq 1/2$ such that $\mE_X (|X_j|^p)^{1/p} \leq \kappa_1 p^{\gamma}$ for all $2 \leq p \leq \kappa_2 \ln(w d)$
($X_j$'s have polynomially growing moments up to the order $\ln d$)
\end{enumerate}
\end{assumption}
The last moments condition of Assumption \ref{as:px}  ensures that
for $n \geq C_1(\ln d)^{\max(2\gamma-1,1)}$, 
\be \label{eq:lm17}
\mE_X\max_{1 \leq j \leq d} \frac{1}{n} \sum_{i=1}^n X_{ij}^2 \leq C_2
\ee
for some constants $C_1=C_1(\kappa_1,\kappa_2,w,\gamma)$ and $C_2=C_2(\kappa_1,\kappa_2,w)$ (Lecu\'e and Mendelson, 2017, proof of Theorem A). Moreover, (\ref{eq:lm17}) might be violated if it holds only up to
the order of $\ln(wd)/\ln\ln(wd)$. In particular, it evidently holds for (scaled) Gaussian and sub-Gaussian $X_j$'s with $\gamma=1/2$ for all moments.   
For simplicity, in what follows we consider
$\gamma=1/2$ corresponding to $n \geq C_1 \ln d$.

\begin{assumption}[{\bf $\mathbf{GSRE(S_0,c_0)}$-condition}] \label{as:GSRE}
For given $S_0$ and $c_0$ assume that 
$$
\kappa(S_0)=\inf_{\Theta_1 \neq \Theta_2 :~ ||\Theta_2-\Theta_1||_1 \leq c_0 \sqrt{S_0} ||\Theta_2-\Theta_1||_2} \frac{||g_{\Theta_2}-g_{\Theta_1}||^2_{L_2(\mP_X)}}{||\Theta_2-\Theta_1||^2_2} >0
$$
\end{assumption}
We refer to Assumpion \ref{as:GSRE} as a generalized strong restricted eigenvalue condition (GSRE), since it can be viewed as a generalization of the strong restricted eigenvalue condition for Lasso in linear regression (e.g., Bellec {\em et al.}, 2018) to the nonlinear regression setup.

We now explore the performance of the Lasso DNN estimator (\ref{eq:lasso}) for regression and classification.

\section{Nonparametric regression} \label{sec:regression}
Consider first a nonparametric regression model with random design:
\begin{equation} \label{eq:regression}
Y_i=g^*(\bX_i)+\epsilon_i,\;\;\;i=1,\ldots,n,
\end{equation}
where $g^*=E(Y|\bX): \mathbb{R}^d \rightarrow \mathbb{R}^1$ is an unknown regression function, $\bX_i \in \mathbb{R}^d \sim \mathbb{P}_X$, $\epsilon_i$'s are zero mean i.i.d. sub-Gaussian random variables independent to $\bX$.
Consider the $l_1$-regularized DNN estimator (\ref{eq:lasso}) with the quadratic loss $\ell(y,g_\Theta(\bx))=(y-g_\Theta(\bx))^2$~:
\begin{equation} \label{eq:lassoreg}
\widehat{\Theta}=\arg\min_{\Theta \in \cF(\ln n,n)}\frac{1}{n} \sum_{i=1}^n (Y_i-g_\Theta(\bx_i))^2 + \lambda ||\Theta||_1
\end{equation}
The goodness of  $\widehat{g}_{\widehat \Theta}$ is measured by the $L_2(\mP_X)$-risk
$$
\mE||\widehat{g}_{\widehat \Theta}-g^*||^2_{L_2(\mP_X)}=\mE\int (\widehat{g}_{\widehat \Theta}(\bx)-g^*(\bx))^2 d\mP_X(\bx)
$$

\subsection{Generalization error bound}
The following theorem established the upper bound for the generalization error of $\widehat{g}_{\widehat \Theta}$:
\begin{theorem}\label{th:general}  Consider the nonparametric regression model (\ref{eq:regression}), where $g^*$ is bounded and $d$ can grow at most polynomially with $n$.
Assume that for any $0 < \epsilon  \leq \epsilon_0$, $g^*$ can be approximated by a sparse DNN $\Theta_\epsilon \in \cF(\ln n,S_\epsilon)$, where 
$S_\epsilon \sim \epsilon^{-\tau} (\ln \epsilon \frac{1}{\epsilon})^r 
$ for some $\tau \geq 0$ and $r \geq 0$, with the $L_2$-approximation error $\epsilon$, that is,
\begin{equation} \label{eq:approx}
\inf_{\Theta \in \cF(\ln n, S_\epsilon)} ||g^*-g_\Theta||_{L_2} \leq \epsilon
\end{equation}

Then, applying the Lasso regularization (\ref{eq:lassoreg}) with $\lambda=C_0 {\sqrt \frac{\ln n}{n}}$  for some $C_0>0$, under Assumption \ref{as:px} and $GSRE(S_\epsilon,10)$-condition, we have
\begin{equation} \label{eq:general}
\mE||\widehat{g}_{\widehat \Theta}-g^*||^2_{L_2(\mP_X)} \leq C \left(\epsilon^{-\tau} \left(\ln \frac{1}{\epsilon}\right)^r\frac{1}{\kappa(S_\epsilon)}~ \frac{\ln n}{n}+\epsilon^2 \right)
\end{equation}
for some $C>0$.
\end{theorem}
The first error term in the RHS of (\ref{eq:general}) bounds the estimation error $\mE||\widehat{g}_{\widehat{\Theta}}-g_{\Theta^*}||^2_{L_2(\mP_X)}$ and is of the order $\frac{S_\epsilon}{\kappa(S_\epsilon)}\frac{\ln n}{n}$, while the second term $\epsilon^2$ is the approximation error $||g^*-g_{\Theta^*}||^2_{L_2}$.

In particular,
taking $\epsilon=n^{-\frac{1}{\tau+2}}$ in (\ref{eq:general}) implies  the immediate corollary:
\begin{corollary} \label{cor:general}
Under all the conditions of Theorem \ref{th:general},
the Lasso DNN estimator (\ref{eq:lassoreg}) with 
$\lambda=C_0 {\sqrt \frac{\ln n}{n}}$ satisfies
$$
E||\widehat{g}_{\widehat \Theta}-g^*||^2_{L_2(\mP_X)} \leq 
C n^{-\frac{2}{2+\tau}} (\ln n)^{r+1} \frac{1}{\kappa(S_0)}
$$
with $S_0 \sim n^\frac{\tau}{\tau+2} (\ln n)^r$.
\end{corollary}

Somewhat similar generalization error bounds were obtained in Ohn and Kim (2022) for a sparse-penalized DNN estimator with the {\em nonconvex} clipped $l_1$-penalty $\lambda \sum_{j=1}^S \max\left(\frac{|\theta_j|}{\gamma},1\right)$, where $\gamma>0$ is the additional tuning parameter. However, nonconvex penalty makes the optimization problem computationally infeasible and, as a practical solution, the authors suggested to replace the original clipped penalty with its convex surrogate.

\subsection{Adaptive minimaxity} \label{subsec:adaptive}
We now demonstrate how the expressivity of DNNs enables to apply the derived general upper bounds for the generalization error of the Lasso DNN estimator to establish its (nearly) minimaxity simultaneously across a wide range of function classes $\cG$. Throughout this section we assume Assumption \ref{as:px} and $GSRE(S_0,10)$-condition with the corresponding $S_0$.

\subsubsection{H\'older and Sobolev spaces of smooth functions}  \label{subsubsec:smooth}
In particular, the approximation condition (\ref{eq:approx}) holds
for H\'older $C^s$
\ignore{
$C^s=\{f: \sum_{\bj: ||\bj||_1 \leq \fs} ||\partial^\bj f||_{L_\infty}+\max_{||\bj||_1=1} \sup_{\bx_1 \neq \bx_2}\frac{|\partial^{\fs}(\bx_1)- \partial^{\fs}(\bx_2)|}{|\bx_1-\bx_2|^{s-\fs}} < \infty\}$
}
and Sobolev $W^{s,p}$
\ignore{
$W^{s,p}=\{f: \sum_{||\bj||_1 \leq s} ||\partial^\bj f||_{L_p} < \infty\}$ for integer $s$
}
classes of smooth functions with regularity $s>0$ for $\tau=d/s$ and $r=1$
(Yarotsky, 2017; Schmidt-Hieber, 2020; G\"uhring {\em et al.}, 2022). Hence, by Corollary \ref{cor:general},  for $g^* \in C^s$ and $g^* \in W^{s,p}$, 
$$
E||\widehat{g}_{\widehat \Theta}-g^*||^2_{L_2(P_X)} =
O\left(n^{-\frac{2s}{2s+d}} \ln^2 n \right) 
$$
which is, up to a log-factor, the minimax quadratic risk in these classes (e.g., Tsybakov, 2009).

Thus, the proposed (nonlinear) $l_1$-regularized DNN estimator $\widehat{g}_{\widehat \Theta}$ is adaptive to the unknown smoothness $s$ and yet is nearly-optimal across the entire range of smooth classes.
On the other hand, there exist various well-known {\em linear}
estimators (e.g., kernel estimators, local polynomial regression, projection estimators) that are also adaptively minimax for $C^s$ and $W^{s,p}$ (Tsybakov, 2009).

We should also mention the remarkable results of Yarotsky (2018) and Yarotsky and Zhevnerchuk (2020) who showed that for the H\'older class $H^{s}$ the approximation error $\epsilon$ in (\ref{eq:approx}) can be achieved with a much smaller number of parameters $S_\epsilon$ with $\tau=d/(2s)$  by very deep DNNs with the number of layers $L_\epsilon \sim \epsilon^{-d/(2s)}$  (compare with $\tau=d/s$ and $L_\epsilon=\ln n \sim \ln \epsilon^{-1}$) although with necessarily discontinuous weight assignments. However, such high expressivity comes at the expense of increasing the estimation error due to accumulating nonlinearities over such a large numbers of layers. The overall quadratic risk evidently cannot be smaller than the minimax risk $n^{-\frac{2s}{2s+d}}$.

\subsubsection{Analytic functions} \label{subsubsec:analytic}
Consider the class of analytic functions - infinitely differentiable functions $g^*$ such that for any $\bx_0 \in \cX$ the Taylor series converges  to $g^*(\bx)$ for any $\bx$ within some neighborhood of $\bx_0$. 

Obshoor {\em et al.} (2022) showed that the condition (\ref{eq:approx}) holds for analytic functions with
$\tau=0$ and $r=d+1$ and, therefore, by Corollary \ref{cor:general},
$$
E||\widehat{g}_{\widehat \Theta}-g^*||^2_{L_2(P_X)} =
O\left(\frac{(\ln n)^{d+2}}{n} \right) 
$$
which is nearly-minimax. 

Note that  similar to H\'older and Sobolev classes, the minimax rates for analytic functions can also be achieved by linear estimators.

\subsubsection{Besov classes} \label{subsec:besov}
The results for H\'older and Sobolev classes from Section \ref{subsubsec:smooth} can be extended to more general Besov classes $B^s_{p,q},\;p,q \geq 1,\;s/d>(1/p-1/2)_+$ of functions. 
The rigorous definition of Besov spaces can be found, e.g. in Meyer (1992). On the intuitive level, (not necessarily integer) $s$ measures the number of function's derivatives, where their existence is required in $L_p$-sense, while $q$ provides a further finer gradation. The Besov spaces include, in particular,
the traditional Sobolev $W^{s,2}$ and H\"older $C^s$ classes of smooth functions
($B^s_{2,2}$ and $B^{s}_{\infty,\infty}$ respectively) but also various spatially inhomogeneous functions for $p<2$
(Meyer, 1992, Donoho and Johnstone, 1998).

Suzuki (2019) showed that $g^* \in B^{s}_{p,q}$ satisfies (\ref{eq:approx}) with
$\tau=d/s$ and $r=1$ and, therefore, by Corollary \ref{cor:general},
$$
E||\widehat{g}_{\widehat \Theta}-g^*||^2_{L_2(P_X)}= 
O\left(n^{-\frac{2s}{2s+d}} \ln^2n \right)
$$
which is the nearly-minimax risk
for $B^{s}_{p,q}$ (Donoho and Johnstone, 1998). Moreover, Donoho and Johnstone (1998) proved  that for $p<2$ linear estimators are only sub-optimal with the best possible risk of the order
$n^{-\frac{s-d/p+d/2}{s-d/p+d}}$. On the other hand, there exist nonlinear estimators based on wavelet thresholding which are minimax across the entire range of Besov spaces including $p<2$.

\subsubsection{Piecewise smooth functions} \label{subsubsec:piecewise}
Imaizumi and Fukumizu (2022) considered a certain class of piecewise smooth functions $\cF^{PS}_{s,\beta,M}$ 
with support partitioned into a (finite)  number of pieces $M$ with $\beta$-smooth boundaries and which are $s$-smooth within each of the pieces (see their paper for the rigorous definition).
They showed that (\ref{eq:approx}) holds for $\cF^{PS}_{s,\beta,M}$ with $\tau=\max\left(\frac{d}{s},\frac{2(d-1)}{\beta}\right)$ and $r=1$. Corollary \ref{cor:general} yields then
$$
E||\widehat{g}_{\widehat \Theta}-g^*||^2_{L_2(P_X)} 
= O\left(\max\left(n^{-\frac{2s}{2s+d}},n^{-\frac{\beta}{\beta+d-1}}\right)\ln^2 n\right)
$$
which is a nearly-minimax risk over $\cF^{PS}_{s,\beta,M}$ (Imaizumi and Fukumizu, 2022). Moreover, Imaizumi and Fukumizu (2022) proved that linear estimators can achieve the risk only of the order $O\left(n^{-\frac{\beta}{2\beta+d-1}}\right)$ which is sub-optimal if $\beta < \frac{2s(d-1)}{d-2s}$. Hence, similar to spatially inhomogeneous functions from Besov classes with $p<2$,
the DNN estimator $\widehat{g}_{\widehat \Theta}$ outperforms linear estimators for piecewise smooth functions with not sufficiently smooth boundaries.

\subsubsection{Composition structured function} \label{subsubsec:composition}
Assume now that the unknown $g^*$ has a certain composition structure. Namely, let $g^*=g^*_L \circ \ldots \circ g^*_1$,
where $g^*_l: [a_l,b_l]^{d_l} \rightarrow [a_{l+1},b_{l+1}]^{d_{l+1}}$ with $t_l \leq d_l$-variable components $g^*_{lj},\;j=1,\ldots, d_l$. Such functions include, in particular, additive, generalized 
additive and single index models (Schmidt-Hieber, 2020). Similar composition structured functions were also considered in Kohler and Langer (2021).
Let $g^*_{lj} \in H^{s_l}$ and define $s_l^*=s_l \prod_{l'=l+1}^L \min(s_{l'},1)$. Schmidt-Hieber (2020)
showed that under mild conditions, such class of composition structured functions satisfies (\ref{eq:approx}) with
$\tau=\max_{1 \leq l \leq L}\frac{t_l}{s_l^*}$ and $r=1$, and, therefore, applying Corollary \ref{cor:general} we have
$$
E||\widehat{g}_{\widehat \Theta}-g^*||^2_{L_2(P_X)}= 
O\left(\max_{1 \leq l \leq L} n^{-\frac{2s_l^*}{2s_l^*+t_l}} \ln^2n \right)
$$
which is nearly-minimax.

The results of Section \ref{subsec:adaptive} are summarized in Table \ref{tab:rateGSREg}. We conclude that the DNN estimator with Lasso regularization (\ref{eq:lassoreg}) is inherently adaptive to the unknown function class and is simultaneously nearly-minimax across the wide range of various smooth and non-smooth functions. 

\begin{table}
\begin{center}

\begin{tabular}{llll} 
\hline
Class  &  $\tau$ & $r$ &  Quadratic risk \\
\hline
H\'older $H^s$, Sobolev $W^{s,r}$  & $\frac{d}{s}$ & 1 & $n^{-\frac{2s}{2s+d}} \ln ^2n$ \\ \\
analytic functions & $0$ & $d+1$ &  $\frac{(\ln n)^{d+2}}{n}$ \\ \\
Besov $B^{s}_{p,q},\;s/d>(1/p-1/2)_+$ & $\frac{d}{s}$ & 1  & $n^{-\frac{2s}{2s+d}} \ln^2 n$ \\ \\
piecewise smooth $\cF^{PS}_{s,\beta,M}$ & $\max\left(\frac{d}{s},\frac{2(d-1)}{\beta}\right)$  & 1 & $\max\left(n^{-\frac{2s}{2s+d}},n^{-\frac{\beta}{\beta+d-1}}\right)\ln^2 n$ \\ \\
composition structured functions & $\max_{1 \leq l \leq L} \frac{t_l}{s_l^*}$ & 1 & $\max_{1 \leq l \leq L} n^{-\frac{2s_l^*}{2s_l^*+t_l}} \ln^2 n$ \\
\hline
\end{tabular}
\caption{Quadratic risks for various function classes.}
\label{tab:rateGSREg}
\end{center}
\end{table}

\section{Classification} \label{sec:classification}
The similar results established for $l_1$-regularized DNN estimator for regression in Section \ref{sec:regression}, also hold for classification. 

\subsection{Binary classification} \label{subsec:binary}
Consider a binary classification setup:
\be \label{eq:classification}
Y|(\bX=\bx) \sim B(1,p^*(\bx)),\;\;\; \bX \sim \mP_X,
\ee
where the logit function $g^*(\bx)=\ln\frac{p^*(\bx)}{1-p^*(\bx)}$ belongs to some class of functions $\cG$ and is bounded. The boundedness of $g(\bx)$ is satisfied, for example, when the conditional probability $p^*(\bx)$ is bounded away from zero and one, i.e. there exists $0 < \delta < \frac{1}{2}$ such that $\delta \leq p^*(\bx) \leq 1-\delta$ for all
$\bx \in \cX$. In particular, it prevents the conditional variance $Var(Y|\bX=\bx)=p^*(\bx)(1-p^*(\bx))$ from being infinitely close to zero.

The quality of a classifier $\eta: \mathbb{R}^d \rightarrow \{0,1\}$ is quantified by the misclassification error $R(\eta)=P(Y \neq \eta(\bX))$ and the optimal Bayes classifier is
$\eta^*(\bx)=I(p^*(\bx) \geq 1/2)$ or, equivalently, $\eta^*(\bx)=I(g^*(\bx) \geq 0)$ with the risk $R(\eta^*)=E_\bX \min(p^*(\bx),1-p^*(\bx))$. Given the data $(\bX_1,Y_1),\ldots,(\bX_n,Y_n)$ one estimates the unknown $g^*$ by some $\widehat{g}$ and designs a plug-in classifier $\widehat{\eta}(\bx)=I\{\widehat{g}(\bx) \geq 0\}$. The goodness of $\widehat{\eta}$ is measured by the misclassification excess risk $\cE(\widehat{\eta},\eta^*)=\mE R(\widehat{\eta})-R(\eta^*)$.

Similar to the regression setup we estimate $g$ by $l_1$-regularized ERM DNN estimator (\ref{eq:lasso}) but w.r.t. the logistic loss $\ell(y,g_\Theta(\bx))=\ln(1+e^{g_\Theta(\bx)})-yg_\Theta(\bx)$ or, in a simpler form,  $\ell(y^*,g_\Theta(\bx))=\ln(1+e^{-y^*g_\Theta(\bx)})$, where
$y^*=2y-1 \in \{-1,1\}$.
Namely,
\begin{equation} \label{eq:lassoclass}
\widehat{\Theta}=\arg\min_{\Theta \in \cF(\ln n,n)}\frac{1}{n}\sum_{i=1}^n \ln\left(1+e^{-Y^*_ig_\Theta(\bX_i)}\right) + \lambda ||\Theta||_1
\end{equation}
Consider the corresponding classifier 
\be \label{eq:classifier}
\widehat{\eta}_{\widehat \Theta}(\bx)=I\{\widehat{g}_{\widehat \Theta}(\bx) \geq 0\}
\ee

Here are the analogs of Theorem \ref{th:general} and Corollary \ref{cor:general} for
classification:
\begin{theorem}\label{th:general1}  Consider the classification model (\ref{eq:classification}), where the logit function $g$ is bounded and
 $d$ can grow at most polynomially with $n$. Assume that for any $0 < \epsilon  \leq \epsilon_0$, there exists a sparse DNN $\Theta_\epsilon \in \cF(\ln n,S_\epsilon)$,  where 
$S_\epsilon \sim \epsilon^{-\tau} (\ln \frac{1}{\epsilon})^r 
$ for some $\tau \geq 0$ and $r \geq 0$, such that
\begin{equation} \label{eq:approx1}
\inf_{\Theta \in \cF(\ln n, S_\epsilon)} ||g^*-g_\Theta||_{L_2} \leq \epsilon
\end{equation}
for any $0 < \epsilon \leq \epsilon_0$.

Then, under Assumption \ref{as:px} and $GSRE(S_\epsilon,10)$-condition, for  the Lasso DNN estimator (\ref{eq:lassoclass}) with $\lambda=C_0 {\sqrt \frac{\ln n}{n}}$  for some $C_0>0$, we have
\begin{equation} \label{eq:general11}
\cE(\widehat{\eta}_{\widehat \Theta}, \eta^*)\leq 
C_1\sqrt{\left(\epsilon^{-\tau} \left(\ln \frac{1}{\epsilon}\right)^r \frac{1}{\kappa(S_\epsilon)}~ \frac{\ln n}{n}+\epsilon^2 \right)}
\end{equation}
and
\begin{equation} \label{eq:general12}
E||\widehat{g}_{\widehat \Theta}-g^*||^2_{L_2(\mP_X)} \leq C_2 \left(\epsilon^{-\tau} \left(\ln \frac{1}{\epsilon}\right)^r \frac{1}{\kappa(S_\epsilon)}~ \frac{\ln n}{n}+\epsilon^2 \right)
\end{equation}
for some $C_1, C_2>0$.
\end{theorem}

Applying Theorem \ref{th:general1} for $\epsilon=n^{-\frac{1}{\tau+2}}$ yields:
\begin{corollary} \label{cor:general1}
Under all the conditions of Theorem \ref{th:general1},
the Lasso DNN estimator (\ref{eq:lassoclass}) with 
$\lambda=C_0 {\sqrt \frac{\ln n}{n}}$ satisfies
\begin{equation} \label{eq:cor11}
\cE(\widehat{\eta}_{\widehat \Theta},\eta^*)\leq 
C_1 n^{-\frac{1}{2+\tau}}(\ln n)^{\frac{r+1}{2}}
\sqrt{\frac{1}{\kappa(S_0)}}
\end{equation}
and
$$
E||\widehat{g}_{\widehat \Theta}-g^*||^2_{L_2(\mP_X)} \leq 
C_2 n^{-\frac{2}{2+\tau}} (\ln n)^{r+1} \frac{1}{\kappa(S_0)}
$$
with $S_0 \sim n^\frac{\tau}{\tau+2} (\ln n)^r$.
\end{corollary}

To give more insight into the upper bound (\ref{eq:cor11}) for the misclassification excess risk, note that the VC-dimension of the set of DNNs $\cF(L_0,S_0)$ is
$$
VC(\cF(L_0,S_0)) \sim  S_0 L_0 \ln S_0 
$$
(Bartlett {\em et al.}, 2019).
Thus, for $L_0=\ln n$ and $S_0 \sim n^\frac{\tau}{\tau+2} (\ln n)^r$ from Corollary \ref{cor:general1}, the bound (\ref{eq:cor11}) is of the order  $\sqrt{\frac{VC(\cF(L_0,S_0))}{n}}$ which is known to be the best possible rate for classifiers from $\cF(L_0,S_0)$ (Devroy {\em et al.}, 1996).

As is done in Section \ref{subsec:adaptive}, we can apply the general bounds in  Corollary \ref{cor:general1} to establish the error rates
of the Lasso-penalized DNN classifier (\ref{eq:classifier}) over various function classes considered there. See Table \ref{tab:ratesclass} for the the resulting misclassification excess risks together with the corresponding $\tau$ and $r$ for the approximation condition (\ref{eq:approx1}) (see Section \ref{sec:regression}).

\begin{table}
\begin{center}

\begin{tabular}{llll} 
\hline
Class  &  $\tau$ & $r$ &  Misclassification excess risk \\
\hline
H\'older $H^s$, Sobolev $W^{s,r}$  & $\frac{d}{s}$ & 1 & $n^{-\frac{s}{2s+d}} \ln n$ \\ \\
analytic functions & $0$ & $d+1$ &  $\sqrt{\frac{(\ln n)^{d+2}}{n}}$ \\ \\
Besov $B^{s}_{p,q},\;s/d>(1/p-1/2)_+$ & $\frac{d}{s}$ & 1  & $n^{-\frac{s}{2s+d}} \ln n$ \\ \\
piecewise smooth $\cF^{PS}_{s,\beta,M}$ & $\max\left(\frac{d}{s},\frac{2(d-1)}{\beta}\right)$  & 1 & $\max\left(n^{-\frac{s}{2s+d}},n^{-\frac{\beta/2}{\beta+d-1}}\right)\ln n$ \\ \\
composition structured functions & $\max_{1 \leq l \leq L} \frac{t_l}{s_l^*}$ & 1 & $\max_{1 \leq l \leq L} n^{-\frac{s_l^*}{2s_l^*+t_l}} \ln n$ \\
\hline
\end{tabular}
\caption{Misclassification excess risks for various function classes.}
\label{tab:ratesclass}
\end{center}
\end{table}

The obtained rates in Table \ref{tab:ratesclass} are nearly-minimax (up to log-factors) for H\'older, Sobolev and Besov classes (Yang, 1999), and we conjecture that this holds true for other considered classes as well. The rates can be improved under the additional low-noise (aka Tsybakov) condition and stronger assumptions on the marginal distribution $\mP_X$ (Audibert and Tsybakov, 2007; Kim {\em et al.}, 2021). However, we leave these issues beyond the scope of the paper.

\subsection{Multiclass classification} 
\label{subsec:multiclass}
The results for binary classification from Section \ref{subsec:binary} can be generalized to the multiclass setup with $K$ classes. In this case,
\be \label{eq:multinom}
Y|(\bX=\bx) \sim Mult(p^*(\bx),\ldots,p^*_K(\bx)),\;\;\;\sum_{k=1}^K p^*_k(\bx)=1,\;\;\; \bX \sim \mP_X,
\ee
where 
\be \label{eq:multinom1}
p^*_k(\bx)=\frac{e^{g^*_k(\bx)}}{\sum_{k'=1}^K e^{g^*_{k'}(\bx)}},\;\;\;\;k=1,\ldots,K
\ee
The model (\ref{eq:multinom})-(\ref{eq:multinom1}) is not identifiable without an extra constraint since shifting each $g^*_k$ by the same $c$ does not affect the probabilities $p^*_k$'s. Without loss of generality, consider the constraint $g^*_K(\bx) \equiv 0$. In this case the last $K$-th class is used as a reference one and
$$
g^*_k(\bx)=\ln \frac{p^*_k(\bx)}{p^*_K(\bx)},\;\;\;k=1,\ldots,K-1
$$
The misclassification risk of a classifier $\eta: \cX \rightarrow \{1,\ldots,K\}$ is $R(\eta)=P(Y \neq \eta(\bX))$ and the Bayes classifier $\eta^*(\bx)=\arg\max_{1 \leq k \leq K}p^*_k(\bx)$ or, equivalenlty, $\eta^*(\bx)=\arg\max_{1 \leq k \leq K} g^*_k(\bx)$, with the oracle risk $R(\eta^*)=1-\mE_X\max_{1 \leq k \leq K} p^*_k(\bX)$.

The DNN architecture used for the multiclass classification is essentially the same as that for regression and binary classification but with $d_L=K-1$ and a vector-valued output $g_\Theta(\bx)=(g_{\Theta 1}(\bx),\ldots,g_{\Theta K-1}(\bx))^T$ in (\ref{eq:dnn}). 

Consider the multinomial logistic loss
$$
\ell(\by,g_\Theta(\bx))=\ln(1+e^{\sum_{k=1}^{K-1}g_{\Theta k}(\bx)})-\sum_{k=1}^{K-1}\xi_k g_{\Theta k}(\bx),
$$
where the indicators $\xi_k=I\{y=k\},\;k=1,\ldots,K-1$, and the corresponding $l_1$-regularized DNN estimator
\be \label{eq:lassomulticlass}
\widehat{\Theta}=\arg\min_{\Theta \in \cF(\ln n,n)}
\frac{1}{n}\sum_{i=1}^n \left(\ln(1+e^{\sum_{k=1}^{K-1}g_{\Theta k}(\bX_i)})-\sum_{k=1}^{K-1}\xi_{ik} g_{\Theta k}(\bX_i)\right) + \lambda ||\Theta||_1
\ee
The resulting classifier $\widehat{\eta}_{\widehat \Theta}(\bx)=\arg\max_{1 \leq k \leq K} \widehat{g}_{\widehat{\Theta}k}(\bx)$, where recall that
$\widehat{g}_{\widehat{\Theta}K}(\bx) \equiv 0$.

We have the following multiclass extension of Theorem \ref{th:general1}:
\begin{theorem}\label{th:general2}  Consider the multiclass classification model (\ref{eq:multinom})-(\ref{eq:multinom1}), where the logit functions $g_k$ are bounded and
 $d$ can grow at most polynomially with $n$. Assume that for any $0 < \epsilon  \leq \epsilon_0$, $g$ there exists a sparse DNN $\Theta_\epsilon \in \cF(\ln n,S_\epsilon)$, where  
$S_\epsilon \sim \epsilon^{-\tau} (\ln\frac{1}{\epsilon})^r 
$ for some $\tau \geq 0$ and $r \geq 0$, such that
\begin{equation} \label{eq:approx2}
\inf_{\Theta \in \cF(\ln n, S_\epsilon)} ||g^*_k-g_{\Theta k}||_{L_2} \leq \epsilon
\end{equation}
for all $1 \leq k \leq K-1$ and any $0 < \epsilon \leq \epsilon_0$.

Then, under Assumption \ref{as:px} and $GSRE(S_\epsilon,10)$-condition, for  the Lasso DNN estimator (\ref{eq:lassomulticlass}) with $\lambda=C_0 \sqrt{(K-1)\frac{\ln n}{n}}$  for some $C_0>0$, we have
$$
\cE(\widehat{\eta}_{\widehat \Theta},\eta^*)\leq 
C_1\sqrt{(K-1) \left(\epsilon^{-\tau} \left(\ln \frac{1}{\epsilon}\right)^r \frac{1}{\kappa(S_\epsilon)}~ \frac{\ln n}{n}+\epsilon^2 \right)}
$$
and
$$
E||\widehat{g}_{\widehat \Theta}-g^*||^2_{L_2(P_X)} \leq C_2 (K-1)\left(\epsilon^{-\tau} \left(\ln \frac{1}{\epsilon}\right)^r \frac{1}{\kappa(S_\epsilon)}~ \frac{\ln n}{n}+\epsilon^2 \right)
$$
for some $C_1, C_2>0$.
\end{theorem}
Taking again $\epsilon=n^{-\frac{1}{\tau+2}}$ leads to the following corollary: 
\begin{corollary} \label{cor:general2}
Under all the conditions of Theorem \ref{th:general2},
the Lasso DNN estimator (\ref{eq:lassomulticlass}) with 
$\lambda=C_0 \sqrt{(K-1)\frac{\ln n}{n}}$ satisfies
$$
\cE(\widehat{\eta}_{\widehat \Theta},\eta^*)\leq 
C_1 n^{-\frac{1}{2+\tau}}(\ln n)^{\frac{r+1}{2}}
\sqrt{\frac{K-1}{\kappa(S_0)}}
$$
and
$$
E||\widehat{g}_{\widehat \Theta}-g^*||^2_{L_2(\mP_X)} \leq 
C_2 n^{-\frac{2}{2+\tau}} (\ln n)^{r+1} \frac{K-1}{\kappa(S_0)}
$$
with $S_0 \sim n^\frac{\tau}{\tau+2} (\ln n)^r$.
\end{corollary}
Similar to the binary case, Corollary \ref{cor:general2} allows one to derive error bounds when $g_k \in \cG$ for various function classes $\cG$ with the corresponding $\tau$.

The misclassification excess risk in Theorem \ref{th:general2} and Corollary \ref{cor:general2}  increases as the square root of the number of classes $\sqrt{K}$, known to be the best possible rate for linear and some other families of classifiers (Danieli {\em et al.}, 2012; Abramovich {\em et al.}, 2021). We conjecture it holds true for multiclass DNN classifiers as well (see also Jin, 2023). On the other hand, based on the results of Jin (2023), it follows that for tree-based classifiers the corresponding price is only of the order $\sqrt{\ln K}$ potentially indicating higher complexity of multiclass classification with DNN classifiers compared to their trees-based counterparts. However, the latter may not possess such expressive power.

\section{Concluding remarks and further challenges} \label{sec:remarks}
We conclude the paper with several remarks and mention some further possible challenges.

\subsection{Various types of DNN sparsity} \label{subsec:sparsity} 
As we have mentioned in Introduction, there are several possible ways to define DNN sparsity. Whereas in this paper we studied connection sparsity (small number of active connections between nodes), one can consider other notions of sparsity, e.g.  {\em node} sparsity (small number of active nodes) and
{\em layer} sparsity (small number of active layers). 

Consider a general regularized DNN estimator
$$
\widehat{\Theta}=\arg\min_{\Theta \in \cF(\ln n,n)}\frac{1}{n} \sum_{i=1}^n \ell\left(Y_i,g_\Theta(\bX_i)\right) + Pen(\Theta)
$$
with a (convex) regularization penalty $Pen(\Theta)$ corresponding to the specific type of sparsity at hand.

The connection sparsity considered in the paper is induced by the Lasso regularized penalty $Pen(\Theta)=\lambda||\Theta||_1$.

To capture node sparsity note that
a $j$-th node at the $l$-layer is not active iff $\sum_{j=1}^{d_{l-1}}W_{l,ij}=0$, i.e. the $i$-th row of $W_l$ is zero. A natural penalty that prompts row-wise sparsity is a {\em group Lasso} 
$$
Pen(\Theta)=\lambda \sum_{l=1}^L ||W_l||_{2,1}=
\lambda \sum_{l=1}^L \sum_{i=1}^{d_l}\sqrt{\sum_{j=1}^{d_{l-1}}W_{l,ij}^2}
$$
(Lederer, 2023).

Similarly, layer sparsity is captured by another group Lasso-type penalty 
$$
Pen(\Theta)=\lambda\sum_{l=1}^L ||W_l||_F,
$$
where $||W_l||_F$ is a Frobenius norm of $W_l$
(Wen {\em et al.}, 2016). 

It is also possible to combine several types of sparsity. Thus, both node and connection sparsity are induced by double row-wise sparsity in matrices $W_l$'s, where their nonzero rows have few nonzero entries. For such type of sparsity consider {\em sparse group Lasso} penalty
$$
Pen(\Theta)=\lambda_1 \sum_{l=1}^L ||W_l||_{2,1} +
\lambda_2 ||\Theta||_1
$$
with the tuning parameters $\lambda_1, \lambda_2>0$.

See Wen {\em et al.} (2016) for other possible notions of sparsity in DNNs.

\subsection{Curse of dimensionality}
The minimax convergence rates over various function classes discussed in Section \ref{subsec:adaptive} suffer from a well-known curse of dimensionality phenomenon -- they rapidly slow down as $d$ increases and require exponentially growing sample sizes for consistent estimation. To handle this problem one needs extra structural assumptions on considered functions, e.g. certain composition structure
(see Section \ref{subsubsec:composition}) or mixed smoothness (Suzuki, 2019). Another reasonable assumption in various applications is that the original high-dimensional data actually resides on a lower dimensional manifold of dimensionality $d_0 \ll d$ (Shaham {\em et al.}, 2018).

To reduce the large number of parameters in a fully-connected DNN (\ref{eq:dnn}) one can consider specific types of sparse architectures.  Thus, in {\em convolutional} neural networks (CNN), $W_l$'s in (\ref{eq:dnn}) are of a special type of Toeplitz matrices. CNNs lead to tremendous success in image processing. Expressivity properties of CNNs are discussed in G\"uhring {\em et al.} (2022, Section 8). While the proposed $l_1$-regularization can be straightforwardly applied to CNNs, the resulting approximation and estimation errors across various function classes should still be investigated.

\subsection{Optimality of DNN estimators}
Although DNNs are minimax-optimal for approximating and estimating functions within classical function spaces such as  H\'older, Sobolev, analytic or Besov,
they are not unique in possessing such properties. The same rates are also achieved by other existing  estimators (see Section \ref{subsec:adaptive}). 
Gaining a better understanding of the types of functions that are specifically suited for DNNs and can explain their excellent performance in various applications, remains an important challenge.

\section*{Acknowledgments}
The work was supported by the Israel Science Foundation (ISF), Grant ISF-1095/22. The author is grateful to Tomer Levy for valuable remarks.

\section*{Appendix}
Throughout the proofs we use various generic positive constants, not necessarily the same each time they appear even within a single equation.

\subsection*{Proofs of Theorems \ref{th:general} and \ref{th:general1}}
We present a joint proof for both
Theorems \ref{th:general} and \ref{th:general1}.

For a given DNN architecture consider a general $l_1$-regularized estimator (\ref{eq:lasso})
\be \label{eq:lassogeneral}
\widehat{\Theta}=\arg\min_{\Theta \in \cF(L,S)}\frac{1}{n} \sum_{i=1}^n \ell\left(Y_i,g_\Theta(\bX_i)\right) + \lambda ||\Theta||_1,
\ee
and the excess risk $
\cE(\widehat{g}_{\widehat \Theta},g^*)=\mE R(\widehat{g}_{\widehat \Theta})-R(g^*)$. 

In particular, for the quadratic loss 
$$
\cE(\widehat{g}_{\widehat \Theta},g^*)=\mE||\widehat{g}_{\widehat \Theta}-g^*||^2_{L_2(P_X)},$$ 
whereas for the logistic loss
$$
\cE(\widehat{g}_{\widehat \Theta},g^*)=
\mE\int \left(\frac{e^{g^*(\bx)}}{1+e^{g^*(\bx)}}(g^*(\bx)-\widehat{g}_{\widehat \Theta}(\bx))+\ln\left(\frac{1+e^{{\widehat g}_{\widehat \Theta}(\bx)}}{1+e^{g^*}(\bx)}\right)\right)dP_X(\bx)
$$

For both losses,
$$
\cE(\widehat{g}_{\widehat \Theta},g^*)=
\left(\mE R(\widehat{g}_{\widehat \Theta})-R(g_{\Theta^*})\right) + \left(R(g_{\Theta^*})-R(g^*)\right),
$$
where $\Theta^*=\arg \min_{\Theta \in \cF(L,S)} R(g)$.

We now go along the lines with Lecu\'e and Mendelson (2018) and Alquier {\em et al.} (2019)  with necessary extensions for DNN estimators to bound the estimation error 
$
\cE(\widehat{g}_{\widehat \Theta},g_{\Theta^*})=  \mE R(\widehat{g}_{\widehat \Theta})-R(g_{\Theta^*}).
$

Define first the following quantities.
Let $\widehat{Rad}(\cF^*(L,S))$ be the empirical Rademacher complexity of $\cF^*(L,S)=\{\Theta \in \cF(L,S): ||\Theta||_1 \leq 1\}$, namely,
\be \label{eq:rademacher}
\widehat{Rad}(\cF^*(L,S))=\mE\left\{\frac{1}{\sqrt n}\sup_{\Theta: ||\Theta||_1 \leq 1 } \sum_{i=1}^n \sigma_i g_\Theta(\bX_i)\Big|\bX_1=\bx_1,\ldots,\bX_n=\bx_n \right\},
\ee
where $\sigma_i,\;i=1,\ldots,n$ are i.i.d.
Rademacher random variables with $P(\sigma_i=1)=P(\sigma_i=-1)=1/2$, and
$$
Rad(\cF^*(L,S))=\mE_X\left\{\widehat{Rad}(\cF^*(L,S))\right\}
$$
be the Rademacher complexity of $\cF^*(L,S)$. Note that the Rademacher complexity is usually defined in the literature with
a normalizing factor $\frac{1}{n}$ in (\ref{eq:rademacher}) instead of $\frac{1}{\sqrt n}$ but we use the latter for convenience following Alquier {\em et al.} (2019).

Define a {\em complexity function}
$$ 
\complexity(\rho) = \sqrt{\frac{C_0 Rad(\cF^*(L,S)) \rho}{2  \sqrt{n}}}
$$
for $\rho>0$, where the exact value of $C_0>0$ is specified in
Alquier {\em et al.} (2019).

Define 
$\mathcal{T}(\rho) = \{\Theta^\prime \in \cF(L,S): ||\Theta^\prime||_1= \rho,~ ||g_{\Theta^\prime}||_{L_2(P_X)} \leq r(2\rho)\}$.
Let $\Gamma_{\Theta^*}(\rho) = \bigcup_{\Theta^\prime : ||\Theta^\prime - \Theta^*||_1 < \frac{\rho}{20}} \partial ||\Theta^\prime||_1$, where the subdifferential $\partial ||\Theta^\prime||_1$ of $||\btheta^\prime||_1$ is $\{{\bf u} \in \mathbb{R}^s: u_i={\rm sign}(\theta^\prime_i),\;\theta^\prime_i \neq 0\;{\rm and}\; u_i \in [-1,1]\;{\rm otherwise}\}$.
The {\em sparsity parameter}  is
   $$
    \sparsity(\rho) =\inf_{\btheta^\prime \in \mathcal{T}(\rho)} \sup_{{\bf u} \in \Gamma_{\Theta^*}(\rho)} {\bf u}^T \btheta^\prime.
    $$
Finally, let $\rho^*$ be any solution of the sparsity inequality
\be \label{eq:sparsity_equtation}
\sparsity(\rho^*) \geq \frac{4}{5}\rho^*
\ee

Applying the results of Lecu\'e and Mendelson (2017) for the quadratic risk and similar results 
of Alquier {\em et al.} (2019) for the logistic risk and bounded $g^*$, for (\ref{eq:lassogeneral}) with $\lambda=C_0 \frac{Rad(\cF^*(L,S))}{\sqrt n}$ we have
	\begin{equation} \label{eq:general_upper}
\cE(\widehat{g}_{\widehat \Theta},g_{\Theta^*}) \leq C_1 \frac{Rad(\cF^*(L,S)) \rho^*}{\sqrt{n}}
       \end{equation}
 and      
   \begin{equation} \label{eq:general_upper1}
    \mE\|\widehat{g}_{\widehat \Theta}-g_{\Theta_*}|^2_{P_X} \leq C_2 \frac{Rad(\cF^*(L,S)) \rho^*}{\sqrt{n}}
    \end{equation}
  for some $C_1, C_2>0$. Evidently (\ref{eq:general_upper}) and (\ref{eq:general_upper1}) coincide for the quadratic loss.

We now find the solution $\rho^*$ of the sparsity inequality (\ref{eq:sparsity_equtation}). 
Let $||\Theta^*||_0=||\btheta^*||_0=S_0$ and $\cJ$ be the corresponding set of indexes of $S_0$ nonzero entries of $\btheta^*$. 
Consider any $\Theta^\prime \in \cF(L,S)$ such that 
$||\Theta^*-\Theta^\prime||_1= \rho^*$ and $||g_{\Theta^*}-g_{\Theta^\prime}||_{L_2(P_X)} \leq r(2\rho^*) =\sqrt{\frac{C_0 Rad(\cF^*(L,S)) \rho^*}{\sqrt n}}$, where $\rho^*$ will be specified later. 
Define $\bu \in \mathbb{R}^S$ such that
$u_j={\rm sign}(\theta^*_j),\;j \in \cJ$ and complete $u_j=\pm 1,\; j \in \cJ^c$ in such a way to ensure 
$$
\sum_{j \in \cJ^c} u_j(\theta_j^\prime-\theta_j^*)=
\sum_{j \in \cJ^c} u_j \theta_j^\prime \geq \sum_{j=S_0+1}^S |\btheta^\prime-\btheta^*|_{(j)}
$$
Evidently, $\bu \in \partial ||\Theta^*||_1 \subseteq \Gamma_{\Theta^*}(\rho^*)$ and
\begin{equation} \label{eq:scalar}
\begin{split}
\bu^T(\btheta^\prime-\btheta^*) & =\sum_{j \in \cJ}
u_j (\theta^\prime_j-\theta^*_j) + \sum_{j \in \cJ^c}
u_j (\theta^\prime_j-\theta^*_j) \geq \sum_{j \in \cJ^c}
u_j (\theta^\prime_j-\theta^*_j) -\sum_{j=1}^{S_0}
|\btheta^\prime-\btheta^*|_{(j)} \\
& \geq ||\btheta^\prime-\btheta^*||_1 - 2\sum_{j=1}^{S_0}|\btheta^\prime-\btheta^*|_{(j)}
\geq \rho^*-2 \sqrt{S_0}||\btheta^\prime-\btheta^*||_2
\end{split}
\end{equation}
Consider two cases.

\noindent
{\em Case 1.} $||\Theta^\prime-\Theta^*||_1 \geq 10 \sqrt{S_0} ||\Theta^\prime-\Theta||_2$. In this case
(\ref{eq:scalar}) straightforwardly yields
$$
\bu^T(\btheta^\prime-\btheta^*) \geq \rho^*-\frac{1}{5} ||\btheta^\prime-\btheta^*||_1 = \frac{4}{5}\rho^*
$$
for any $\rho^*>0$.

\vspace{.5cm}
\noindent
{\em Case 2.} $||\Theta^\prime-\Theta^*||_1 < 10 \sqrt{S_0} ||\Theta^\prime-\Theta||_2$. By the condition $GSRE(S_0,10)$ we have
$$
\bu^T(\btheta^\prime-\btheta^*) \geq \rho^*-2\sqrt{S_0} \frac{||g_{\Theta^\prime}-g_{\Theta^*}||_{L_2(P_X)}}{\sqrt{\kappa(S_0)}} \geq
\rho^* - C\sqrt{S_0} \sqrt{\frac{Rad(\cF^*(L,S))\rho^*}{\sqrt n}} \frac{1}{\sqrt{\kappa(S_0)}}
$$
Take 
\be \label{eq:rho*}
\rho^*=c_0 \frac{S_0 Rad(\cF^*(L,S))}{\kappa(S_0)\sqrt{n}},
\ee
where we can always find $c_0>0$ to satisfy
$$
C\sqrt{S_0} \sqrt{\frac{Rad(\cF^*(L,S))\rho^*}{\sqrt n}} \frac{1}{\sqrt{\kappa(S_0)}} = \frac{1}{5}\rho^*
$$

Thus, $\rho^*$ in (\ref{eq:rho*}) is the solution of the sparsity inequality (\ref{eq:sparsity_equtation}) for both cases, and from
(\ref{eq:general_upper})-(\ref{eq:general_upper1})  we have
$$
\cE(\widehat{g}_{\widehat \Theta},g_{\Theta^*}) \leq C_1 \frac{Rad(\cF^*(L,S))^2 S_0}{\kappa(S_0) n}
$$
 and
$$
    \mE\|\widehat{g}_{\widehat \Theta}-g_{\Theta_*}|^2_{\mP_X} \leq C_2 \frac{Rad(\cF^*(L,S))^2 S_0}{\kappa(S_0) n}
$$

By Theorem 3.2 of Golowich {\em et al.} (2020), 
$$
\widehat{Rad}(\cF^*(L,S)) \leq \sqrt{\frac{2(L
+1+\ln d)}{n}}~ \sqrt{\max_{1 \leq j \leq d} \sum_{i=1}^n x^2_{ij}} 
$$
and, therefore, under Assumption \ref{as:px}, the Rademacher complexity 
$$
Rad(\cF^*(L,S)) \leq C \sqrt{L+\ln d}
$$
Hence, with $\lambda=C_0 \sqrt{\frac{L+\ln d}{n}}$,
\be \label{eq:stocherror}
\cE(\widehat{g}_{\widehat \Theta},g_{\Theta^*}) \leq C_1 \frac{(L+\ln d)S_0}{\kappa(S_0) n}
\ee
and, similarly,
\be\label{eq:stocherror1}
\mE\|\widehat{g}_{\widehat \Theta}-g_{\Theta_*}|^2_{\mP_X} \leq C_2 \frac{(L+\ln d)S_0}{\kappa(S_0) n}
\ee

By substituting $L=\ln n,\; S_0=S_\epsilon=\epsilon^{-\tau} (\ln \frac{1}{\epsilon})^r$ from Theorems \ref{th:general}-\ref{th:general1} into (\ref{eq:stocherror}) and (\ref{eq:stocherror1}) and observing that for $d$ growing at most polynomially with $n$, $\ln n + \ln d \sim \ln n$, we obtain the bounds for the estimation errors. 

The approximation errors follow directly from the conditions (\ref{eq:approx}) and (\ref{eq:approx1}) of the theorems. Indeed, for an almost surely bounded marginal density $p_X$ (Assumption \ref{as:px}) and the
the quadratic loss,   (\ref{eq:approx}) yields
$$
R(g_{\Theta^*})-R(g^*)=||g_{\Theta^*}-g^*||^2_{L_2(P_X)} \leq C ||g_{\Theta^*}-g^*||^2_{L_2}
\leq \epsilon^2
$$

Similarly, for the logistic loss under the condition (\ref{eq:approx1}), by Taylor expansion (see also Lemma 1 of Abramovich and Grinshtein, 2016) one has
$$
R(g_{\Theta^*})-R(g^*)  \leq C ||g_{\Theta^*}-g^*||^2_{L_2(P_X)} \leq C ||g_{\Theta^*}-g^*||^2_{L_2}
\leq \epsilon^2
$$  
Finally, to establish (\ref{eq:general11}), we exploit
the well-known relations between the misclassification excess risk and the logistic risk:
$\cE({\widehat \eta}_{\widehat \Theta},\eta^*) \leq 
\sqrt{2\cE(\widehat{g}_{\widehat \Theta},g^*)}$ (Bartlett {\em et al.}, 2006).

\subsection*{Proof of Theorem \ref{th:general2}}
The proof of Theorem \ref{th:general2} for multiclass case essentially replicates the proof of Theorem \ref{th:general1} for binary classification incorporating the extended definition of the empirical Rademacher complexity  (\ref{eq:rademacher}) for a vector-valued function:
$$
\widehat{Rad}(\cF^*(L,S))=\mE\left\{\frac{1}{\sqrt n}\sup_{\Theta: ||\Theta||_1 \leq 1 } \sum_{i=1}^n \sum_{k=1}^{K-1} \sigma_{ik} g_{\Theta k}(\bX_i)\Big|\bX_1=\bx_1,\ldots,\bX_n=\bx_n \right\},
$$
where $\sigma_{ik},\;i=1,\ldots,n,\;k=1,\ldots,K-1$ are i.i.d.
Rademacher random variables with $P(\sigma_{ik}=1)=P(\sigma_{ik}=-1)=1/2$. 

By repeating the proof of Theorem 3.2 of Golowich {\em et al.} (2020) with corresponding evident modifications, one can verify that
$$
\widehat{Rad}(\cF^*(L,S))= \leq  \sqrt{\frac{2(L+(K-1)+
(K-1)\ln d)}{n}}~ \sqrt{\max_{1 \leq j \leq d} \sum_{i=1}^n x^2_{ij}} 
$$
and, therefore, under Assumption \ref{as:px},
\be \label{eq:rademachermulti}
Rad(\cF^*(L,S)) \leq C \sqrt{L+(K-1)\ln d}
\ee

Substituting the above bound (\ref{eq:rademachermulti}) for $Rad(\cF^*(L,S))$
into (\ref{eq:rho*}) and subsequently into the general estimation errors bounds (\ref{eq:general_upper}) and (\ref{eq:general_upper1})
with $L=\ln n$ and $S_0=S_\epsilon=\epsilon^{-\tau}(\ln \frac{1}{\epsilon})$ for $\lambda=C_0 \sqrt{\frac{L+(K-1)\ln d}{n}}$ and $d$ growing at most polynomially with $n$, yields 
$$
\cE(\widehat{g}_{\widehat \Theta},g_{\Theta^*}) \leq C_1(K-1) \epsilon^{-\tau} \left(\ln \frac{1}{\epsilon}\right)^r \frac{1}{\kappa(S_\epsilon)}~ \frac{\ln n}{n}
$$
and
$$
\mE\|\widehat{g}_{\widehat \Theta}-g_{\Theta_*}|^2_{\mP_X} \leq C_2 (K-1)\epsilon^{-\tau} \left(\ln \frac{1}{\epsilon}\right)^r \frac{1}{\kappa(S_\epsilon)}~ \frac{\ln n}{n}
$$
Adding the approximation error $\epsilon^2$ completes the error bounds for $\cE(\widehat{g}_{\widehat \Theta},g^*)$ and
$\mE\|\widehat{g}_{\widehat \Theta}-g^*\|^2_{\mP_X}$. Finally,
similar to the binary classification,
$\cE({\widehat \eta}_{\widehat \Theta},\eta^*) \leq 
\sqrt{2\cE(\widehat{g}_{\widehat \Theta},g^*)}$
(Pires and Szepesvari, 2016; Abramovich {\em et al.}, 2021).

\end{document}